# Direct solution of the inverse optimization problem of load sharing between muscles

Adam Siemieński

Biomechanics Laboratory, University School of Physical Education, Wrocław, Poland

**Introduction**: Muscles crossing a joint usually outnumber its degrees of freedom, making it impossible to evaluate their forces unambiguously. This redundancy is typically handled by assuming that loads are shared between muscles according to an optimization law, which leads to constrained minimization of an objective function [3]. It is not clear, however, how to express such function in terms of muscle forces. This could be answered by inverse optimization, an approach realized so far only indirectly [1, 2, 5].

**Objective**: To propose a direct method of constructing an objective function whose minimization results in a given load sharing pattern.

**Methods**: The construction is based on the observation that formulating a linearly constrained load sharing optimization problem by using Lagrange multipliers leads to Schröder's functional equation [4] for the derivative of the objective function. Solving the inverse optimization problem is therefore equivalent to solving this equation, and then integrating its solution.

**Results and conclusions**: For a linearly constrained minimization problem

$$\sum_i r_i F_i = M \qquad K = \sum_i J(F_i) \to Minimum$$

the associated inverse optimization problem consists in finding $K$ based on known force-force functions. It can be shown that if such a function is known for at least one pair of muscles, i.e. a monotone function $h$ is known, $F_j = h(F_k)$, $h(0)=0$, $0 < h'(0) < 1$, then $J'$ must satisfy Schröder's functional equation, $J'(h(x)) = (r_j/r_k) J'(x)$, for which Koenigs's solution [4] is defined by

$$J'(x) = Const \cdot \left[ \lim_{n \to \infty} \frac{h_n(x)}{h'(0)^n} \right]^p$$

with $h_n$ denoting $n$-th iteration of $h$, and $p$ depending on the ratio $r_j/r_k$. A unique up to one multiplicative and one additive constant objective function $K$ can thus be constructed for virtually any differentiable monotone load sharing pattern $h$. Moreover, within this framework, load sharing patterns for other pairs of muscles crossing the same joint turn out to be fractional iterates of $h$.

In a special case of linear $h$ the objective function proposed in [3] is reproduced. The above limit can also be evaluated exactly for some nonlinear functions $h$, but for a general force-force function, e.g. one obtained experimentally, an approximate formula must be sought instead by stopping iteration at some $n$. Practically, due to quick convergence, stopping at $n=10$ usually results in $J'$ indistinguishable from the limit.